\documentclass[a4paper]{amsart}

\usepackage{amssymb,amsfonts,amsmath,latexsym}
\usepackage{color,graphics,epsfig,float}

\title[A matrix interpolation between classical and free max
  operations]{A matrix interpolation between classical and free max
  operations.\\
 I. The univariate case}

\author[F.~Benaych-Georges]{Florent Benaych-Georges}
\address{Florent Benaych-Georges: LPMA, UPMC Univ Paris 6, Case courier 188, 4, Place Jussieu, 75252 Paris Cedex 05, France\\ and CMAP, \'Ecole Polytechnique, route de Saclay, 91128 Palaiseau Cedex, France} \email{florent.benaych@gmail.com}

\author[T.~Cabanal-Duvillard]{Thierry Cabanal-Duvillard}
\address{Thierry Cabanal-Duvillard: MAP 5, UMR CNRS 8145 - Universit\'e Paris Descartes, 45 rue des Saints-P\`eres 75270 Paris cedex~6, France} \email{Thierry.Cabanal-Duvillard@mi.parisdescartes.fr}

\date{\today}
\newcommand{\kinf}{\underset{k\to\infty}{\longrightarrow}}
\newcommand{\dinf}{\underset{N\to\infty}{\longrightarrow}}

\newcommand{\be}{\begin{equation}}
\newcommand{\ee}{\end{equation}}
\newcommand{\z}{\mathbb{Z}}

\newcommand{\vfi}{\varphi}
\newcommand{\ste}{\, ;\, }

\newcommand{\f}{\frac}
\newcommand{\ff}{\frac{1}}
\newcommand{\st}{such that }
\newcommand{\la}{\lambda}

\newcommand{\pro}{probability }
\newcommand{\eps}{\varepsilon}
\newcommand{\maxd}{\operatorname{max}_N}
\newcommand{\maxN}{\operatorname{max}_N}
\newcommand{\ud}{\operatorname{d}\!}

\renewcommand{\P}{\mathbb{P}}

\newcommand{\R}{\mathbb{R}}

\newcommand{\E}{\mathbb{E}}
\newcommand{\fleche}{\longrightarrow}
\def\be{\begin{equation}}
\def\ee{\end{equation}}
\def\bea{\begin{eqnarray}}
\def\eea{\end{eqnarray}}
\def\beas{\begin{eqnarray*}}
\def\eeas{\end{eqnarray*}}
\newtheorem{thm}{Theorem}[section]
\newtheorem{prop}[thm]{Proposition}

\newtheorem{lem}[thm]{Lemma}
\newtheorem{cor}[thm]{Corollary}

\newtheorem{rmq}[thm]{Remark}
\newenvironment{pr}{\noindent {\bf Proof. }}{\ \ \ $\square$}

\long\def\symbolfootnote[#1]#2{\begingroup
\def\thefootnote{\fnsymbol{footnote}}\footnote[#1]{#2}\endgroup} 

\newcommand{\im}{\mbox{Im}}

\begin{document}
\maketitle
\symbolfootnote[0]{{\it MSC 2000 subject classifications.}  15A52, 
46L54} 

\symbolfootnote[0]{{\it Key words.} random matrices, free probability, max-stable laws}

\begin{abstract}Recently, Ben Arous and Voiculescu considered taking the maximum of two free random variables and brought to light a
  deep analogy with the operation of taking the maximum of two independent
  random variables. We present here a new insight on this analogy: its concrete realization based on
  random matrices giving an interpolation between classical and free
  settings.
\end{abstract}

\tableofcontents

\section{Introduction}

The free probability theory  has been a very active field in mathematics over
the last two decades, constructed in a deep analogy with classical
probability theory. Nowadays, there is an unofficial dictionary of concepts
in  both theories:  many fundamental notions or results of classical
probability theory, such as    Law of Large Numbers, Central Limit Theorem,
Gaussian distribution, convolution, cumulants, infinite divisibility have their
precise counterpart in free probability theory.  

Recently,  Ben
Arous and Voiculescu \cite{BAV} have added a new item in this dictionary, bringing to
light the fact that the operation called the {\it classical upper extremal convolution}, which
associates the distribution  of the supremum $X\vee Y$ to the
distributions of two independent random variables $X$
and $Y$, has an analogue in
free probability theory, called the {\it free upper extremal convolution}, which
associates the distribution  of the supremum $X\vee Y$ to the
distributions of two free random variables $X$
and $Y$. From the point of view of  cumulative distribution functions\footnote{The {\it cumulative distribution function} of a random
  variable $X$ with distribution $\mu$ is the function, denoted either by
  $F_X$ or $F_\mu$, is defined by $F_\mu(x)=\mu((-\infty, x])$, for all $x\in \R$.}, both of these operations have a concise interpretation: the classical (resp. free) upper extremal convolution of two \pro measures with distribution functions $F,G$ is the \pro measure with distribution function $FG$ (resp. $\max(0, F+G-1)$). 
The purpose of this paper is to 
construct random matrix models providing
a better understanding of the free upper extremal convolution and a new insight on  the relation between classical and free upper extremal convolutions.

Let us first notice the following fact (which is an immediate consequence of remark \ref{08.03.11-RMQ}):
for $\mu,\nu$ \pro measures on the real line, for $(X_i)_{i\geq 1}$ and
$(Y_i)_{i\geq 1}$ two independent families of independent identically
distributed random variables, with respective distributions   $\mu,\nu$,
if for each $N\geq 1$, one denotes by $Z_{N,1}\geq Z_{N,2}\geq\cdots\geq
Z_{N,N}$ the $N$ largest elements of the multiset $(X_1,\ldots, X_N, Y_1,
\ldots, Y_N)$, then the empirical \pro measure
$\ff{N}\sum_{i=1}^N\delta_{Z_{N,i}}$ converges almost surely to the free
upper extremal convolution of $\mu$ and $\nu$ as $N$ tends to infinity.
This result may appear as a coincidence. But, if we consider $(X_1,\ldots,X_N)$ and $(Y_1,\ldots,Y_N)$ as
the eigenvalues of two independent random matrices $X^{(N)}$ and $Y^{(N)}$ which are invariant in
law under conjugation by any unitary matrix, then
$(Z_{N,1},\ldots,Z_{N,N})$ are the eigenvalues of the supremum\footnote{The {\it supremum} of two Hermitian matrices is defined in section \ref{pj.18.10.08}.} of $X^{(N)}$
and $Y^{(N)}$ with respect to the spectral order introduced by Olson 
\cite{Ol} and used by Ben Arous and Voiculescu in their paper. Since
$X^{(N)}$ and $Y^{(N)}$ are known to be asymptotically free, this gives
actually a first
interpretation of the free upper extremal convolution.

We would like to give 
a deeper insight on the analogy between  the classical and free 
settings. There is a morphism ${\Lambda}^\vee$ from the set of \pro measures endowed with the
classical upper convolution to the same set endowed with the free upper
convolution. It maps any \pro measure with distribution function $F$ to the
one with distribution function $\max (0, 1+\log F)$.
Moreover, Ben Arous and Voiculescu have established that ${\Lambda}^\vee$ provides a remarkable correspondence between classical
max-stable laws and free max-stable laws, which preserves  the domains of
attraction. 
In this paper, we shall give a concrete realization of this morphism  via a
random matrix model.  

More precisely,
for each \pro measure $\mu$ on the real line and each positive integer $N$, we shall define a law ${\Lambda}^\vee_N(\mu)$ on the set of $N$ by $N$ Hermitian matrices such that:
\begin{itemize}
\item[a)] a ${\Lambda}^\vee_N(\mu)$-distributed random matrix is invariant 
under conjugation by any unitary matrix,
\item[b)] for any pair $\mu_1,\mu_2$ of \pro measures on the real line and any pair $M_1,M_2$ of independent random matrices distributed respectively with respect to ${\Lambda}^\vee_N(\mu_1)$ and ${\Lambda}^\vee_N(\mu_2)$, the law of the supremum $M_1\vee M_2$ of $M_1$ and $M_2$ is ${\Lambda}^\vee_N(\mu)$, where $\mu$ is the classical upper extremal convolution of $\mu_1$ and $\mu_2$,
\item[c)] if for all $N$, $M_N$ is a random matrix with law
  ${\Lambda}^\vee_N(\mu)$, then  the empirical spectral law of $M_N$ tends
  almost surely to ${\Lambda}^\vee(\mu)$ as $N$ tends to infinity.\end{itemize}

 Such a model produces a new clue
with regards to
the relevance of the analogy between the two max-operations, and it may be used
to provide more intuitive proofs of some of Ben Arous and Voiculescu's results.
For instance, we have already mentioned that unitarily-invariant random
matrices behave asymptotically as free random variables. Therefore
one can   re-prove immediately, as a consequence of a), b) and c), the fact that 
${\Lambda}^\vee$ is a morphism between the classical and the free upper extremal convolutions (and thus the formula $F_{X\vee Y}=\max(0, F_X+F_Y-1)$ for free random variables $X,Y$). We shall see later that Ben Arous and Voiculescu's  results about max-stable laws and domains of attraction can also be proved with our random matrix model.

In a forthcoming paper, we will show how this approach is also appropriate  to
study the multivariate counterpart of the free upper
extremal convolution.

The paper is organized as follows. First, we shall focus on the upper extremal
convolution for classical, free and Hermitian random variables. Then, we shall introduce our
models, via a strategy which is close to the one of our previous papers \cite{FBG-05} and \cite{CD-05}, and set the main theorems. The last part will be dedicated
to the proofs.

\section{The upper extremal convolutions}

\subsection{For independent real random variables}

Let $X$ and $Y$ be real independent random variables, with cumulative distribution
functions $F_X$ and $F_Y$, and let us denote $X\vee Y=\max(X,Y)$. The cumulative distribution
function of $X\vee Y$
is then equal to $F_{X\vee Y}=F_XF_Y$. 

\subsection{For free random variables}\label{hood.19.10.08}

Let $(M,{\tau})$ be a tracial $W^*$-probability space, that is a von
Neumann algebra $M$ endowed with an ultraweakly continuous faithful
trace-state ${\tau}$. The {\it maximum
$X\vee Y$
of two self-adjoints elements} $X$ and $Y$ in $M$ is
defined with respect to the spectral order introduced by Olson in 1971
(cf  \cite{Ol}; see also \cite{An,BAV}):\\
$\bullet$ If $p$ and $q$ are two self-adjoint projectors, then $q\vee q$ is the
  Hermitian projector on $\overline{\Im(p)+\Im(q)}$.\\
$\bullet$ If $a$ and $b$ are two self-adjoint elements, with resolutions of
  identity $E_a$ and $E_b$, then $a\vee b$ is the self-adjoint element with resolution of
  identity given by $E_a\vee E_b$: $a\vee b$ is the only self-adjoint element $h$ \st for all real number $t$, $$\chi_{(t,\infty)}(h)=\chi_{(t,\infty)}(a)\vee\chi_{(t,\infty)}(b),$$where $\chi_{(t,\infty)}$ denotes the indicator function of the interval $(t,\infty)$ and is applied   {\it via} the functional calculus.

 In 2006, in \cite{BAV}, Ben Arous and Voiculescu have
determined the cumulative distribution
function of $X\vee Y$ whenever $X,Y$ are free:
\be\label{15.06.08.3}
\forall t\in\R,\ \ F_{X\vee Y}(t)=\max(0,F_X(t)+F_Y(t)-1)=:F_X\Box\hskip-9.3pt\vee F_Y(t).
\ee
It is easy to note that the function ${\Lambda}^\vee(u)=\max(0,1+\log u)$ is
a kind of  morphism between classical and free upper extremal
convolutions:
$$
{\Lambda}^\vee(F_XF_Y)={\Lambda}^\vee(F_X)\Box\hskip-9.3pt\vee{\Lambda}^\vee(F_Y).
$$

{\bf Notations:} 
$\bullet$ In order to avoid the use of too many notations, we
shall define  the binary operation $\displaystyle\Box\hskip-7pt\vee$ on the
set of \pro measures on the real line with the same symbol as the operation
$\displaystyle\Box\hskip-7pt\vee$ on the set of cumulative distribution
functions. For all pair $\mu, \nu$ of \pro measures on the real line,
$\mu\Box\hskip-9.3pt\vee\nu$ is defined by:
$$
F_{\mu\displaystyle\Box\hskip-9.3pt\vee\nu}=F_\mu \Box\hskip-9.3pt\vee
F_\nu\textrm{ (}=\max(0,F_\mu(t)+F_\nu(t)-1)\textrm{ ).}
$$ 
$\bullet$ We shall also define   the operator ${\Lambda}^\vee$ on the set of \pro
measures  on the real line with the same symbol as the function
${\Lambda}^\vee(x)=\max(0,1+\log x)$ on $[0,+\infty)$: for all   \pro
measure $\mu$ on the real line, ${\Lambda}^\vee(\mu)$ is the \pro measure
with cumulative distribution function:
$$
F_{{\Lambda}^\vee(\mu)}(t)={\Lambda}^\vee(F_\mu(t))\textrm{
  (}=\max(0,1+\log F_\mu(t))\textrm{ ).}
$$

\begin{rmq}\label{08.03.11-RMQ}
The free upper extremal convolution also appears in a much simpler situation
than the maximum of two free operators. Let
$x_1,\ldots,x_N,y_1,\ldots,y_N\in\R$,
  and $z_1\geq\ldots\geq z_N$ 
the $N$ largest elements of the multiset
$(x_1,\ldots,x_N,y_1,\ldots,y_N)$. Then $\frac{1}{N}\sum_{i=1}^N{\delta}_{z_i}$ is the free upper extremal convolution of
$\frac{1}{N}\sum_{i=1}^N{\delta}_{x_i}$ and $\frac{1}{N}\sum_{i=1}^N{\delta}_{y_i}$.
\end{rmq}

\subsection{For independent Hermitian random matrices}\label{pj.18.10.08}

Let $N$ be a positive integer. The maximum operation introduced by Olson is also defined on the space of
$N\times N$ Hermitian matrices, as a particular case of the definition given in section \ref{hood.19.10.08}.

Concretely, for $A,B$ Hermitian matrices with the same size, $A\vee B$ can be defined in the following way.
Let $\{\la_1>\la_2>\cdots >\la_p\}$ be the union of the spectrums of $A$
and of $B$. For  $i=0,\ldots, p$, let us define $E_i=\sum_{j=1}^i(\ker(A-\la_jI)+
\ker(B-\la_jI))$. Then $A\vee B$ is the Hermitian matrix with eigenvalues the $\la_j$'s \st $E_j\neq E_{j-1}$, with associated eigenspace $E_j\cap(E_{j-1}^\bot)$. In other words, if, for any subspace $F$, $P_F$ designs the orthogonal projector onto $F$, we have $$A\vee B=\la_1P_{E_1}+\la_2P_{E_2\cap(E_1^\bot)}+\la_3P_{E_3\cap(E_2^\bot)}+\cdots+\la_p P_{E_p\cap(E_{p-1}^\bot)}.$$

The following proposition resumes the characteristics of the operation $\vee$ on unitarily invariant random matrices. Recall that the {\it empirical spectral law} of a matrix is the uniform measure on the multiset of its eigenvalues counted by multiplicity.

\begin{prop}\label{21.03.08.propo1}If $A$ and $B$ are two $N\times N$ Hermitian
random independent matrices, whose laws are invariant by conjugation by any
unitary matrix, then
\begin{itemize}
\item $A\vee B$ is an Hermitian
random  matrix, whose law is invariant under the conjugation by any
unitary matrix;
\item the $ N$ eigenvalues of $A\vee B$ are the $ N$ largest of the $2 N$
  eigenvalues of $A$ and $B$ (counted with multiplicity); 
  \item the cumulative distribution functions $F_A,F_B$ and $F_{A\vee B}$
    of the respective empirical spectral laws of $A,B$ and $A\vee B$ are almost surely linked by the relation  \be\label{max.mat.15.06.08} F_{A\vee B}=\max(0, F_A+F_B-1).\ee
\end{itemize}
\end{prop}

\begin{rmq}It is noteworthy that formulas \eqref{15.06.08.3} and
  \eqref{max.mat.15.06.08} are identical, showing
  that the matricial and free upper extremal convolutions are somehow the
  same operation. As free random variables can be approximated by
  independent unitarily invariant random matrices,  the rather trivial result \eqref{max.mat.15.06.08}
  provides an intuitive proof of formula \eqref{15.06.08.3}. 
\end{rmq}

\section{The matricial interpolation}
The aim of this section is to introduce random matrices giving quite a natural
interpretation of the mapping ${\Lambda}^\vee$. We simply follow the same
approach as we did for the Bercovici-Pata bijection in \cite{FBG-05} and
\cite{CD-05}.

Let $F_\mu$ denote its  cumulative
distribution function for $\mu$ \pro measure on the real
line.  
 For all integer $k\geq 1$, $F_\mu^\ff{k}$ is the cumulative distribution
 function of a \pro measure on the real line, which will be denoted by
 $\mu^\ff{k}$. It is clear that if $X_1,\ldots, X_k$ are independent random
 variables distributed according to $\mu^\ff{k}$, $\max \{X_1,\ldots, X_k\}
 $ should be distributed according to $\mu$.

Let us consider $N,k\geq 1$. Let $X_1,\ldots,X_ N$ be independent random variables distributed according to  $\mu^\ff{k}$ and $U$ an unitary, Haar distributed,
$N\times N$
random matrix independent of $X_1,\ldots,X_ N$. Define the random matrix $M$
by
\be\label{18.06.08.1}
M=U\left(
  \begin{array}{ccc}
    X_1&&\\
&\ddots&\\
&&X_N
  \end{array}
\right)U^*
\ee
Let $M_1,\ldots,M_k$ be independent replicas of $M$. 

\begin{thm}\label{21.03.08.thm}As $k$ tends to infinity, the distribution of $M_1\vee\cdots \vee M_k$ converges to a \pro measure $\Lambda_N^\vee(\mu)$ on the space of $N\times N$ Hermitian matrices. This measure is invariant under the action of the unitary group by conjugation and in the case where   $\mu$ has no atom,  under $\Lambda_N^\vee(\mu)$, the joint distribution of the ranked eigenvalues  is absolutely continuous with respect to  $\mu^{\otimes N}$, with density
\begin{equation}\label{density.21.03.08}
N^N
1_{t_1\geq\cdots\geq
  t_N}\prod_{i=1} ^N\frac{F_\mu(t_N)}{F_\mu(t_i)}
\end{equation}at any point $(t_1,\ldots, t_N)\in \R^ N$ \st for all
$i=1,\ldots,  N$, $F_\mu(t_i)\neq 0$ (the density being set to zero anywhere else).

Moreover, 
 if $\nu$ is another \pro measure on the real line and $A,B$ are independent random matrices respectively distributed according to $\Lambda_N^\vee(\mu)$, $\Lambda_N^\vee(\nu)$, then $A\vee B$ is distributed according to $\Lambda_N^\vee(\rho)$, where $\rho$ is the \pro measure on the real line \st $F_\rho=F_\mu F_\nu$.
\end{thm}

\begin{cor}\label{ultima}
  With the preceding notations, we have

  (i) $\Lambda_N^\vee(\mu)$ is max-infinitely divisible,
  
  (ii) if ${\mu}$ is max-stable,  then so is $\Lambda_N^\vee(\mu)$.
  
\end{cor}

In the following theorem, we are interested in the limit of the spectral
measure of  a ${\Lambda}_N^\vee(\mu)$-distributed random matrix, when its
dimension $ N$ tends to infinity.

\begin{thm}\label{bachelorette.17.04.08}
 Let $\mu$ be a \pro measure on the real line. For all $N\geq 1$, let  $M_ N$ a random
 matrix with law ${\Lambda}_N^\vee(\mu)$,
 ${\lambda}_1,\ldots,{\lambda}_ N$ its eigenvalues and
 $\hat\mu_N=\frac{1}{N}\sum_{i=1} ^N{\delta}_{{\lambda}_i}$ its empirical
 spectral law. Then, when $ N$ goes to infinity, $\hat\mu_ N$ converges
 weakly almost surely to  ${\Lambda}^\vee(\mu)$. 
\end{thm}

The previous theorem 
may be used to
provide new, more intuitive proofs of some of the Ben Arous and Voiculescu's results in
\cite{BAV}. 

A first example
has been given
in the introduction: it is a new 
derivation of the formula \eqref{15.06.08.3} defining
the free upper extremal convolution in terms of cumulative distribution functions. 


For second example,
the fact that ${\Lambda}^\vee$ maps any classical max-stable law to a free one is a direct consequence of (ii) of Corollary \ref{ultima} and of formulas \eqref{15.06.08.3} and
  \eqref{max.mat.15.06.08}. 
  
Finally, let us  consider the preservation
 of the domains of attraction. 
 Ben Arous and Voiculescu have proved that for any classical max-stable law $\mu$, $\Lambda^\vee (\mu)$ 
is freely max-stable and that for any cumulative distribution function $F$,
for any sequences $a_k>0, b_k\in \R$, we have
\be\label{sht.18.10.08}F(a_k\cdot+b_k)^k\kinf F_\mu(\cdot)\Longrightarrow
F(a_k\cdot+b_k)^{\Box\hskip-5.6pt\vee k}\kinf \Lambda^\vee
(F_\mu(\cdot))\ee(to state the reciprocal implication, one needs to choose
a right inverse to $\Lambda^\vee$, which we shall not do here). This result
can be generalized as follows: for all sequence $(F_k)$ of cumulative
distribution functions, if $F$ is a cumulative distribution function, we
have  \be\label{sht.18.10.08.1} F_k(\cdot)^k\kinf F(\cdot)\Longrightarrow
F_k(\cdot)^{\Box\hskip-5.6pt\vee k}\kinf \Lambda^\vee (F(\cdot))\ee(the
proof of \eqref{sht.18.10.08.1} is given in the last section of the
paper). 

Let us now explain how our random matrix model gives an heuristic 
explanation of 
\eqref{sht.18.10.08.1}. 
For sake of simplicity, let $F_k$ be equal to $F^{\frac{1}{k}}$, with $F$ the
cumulative distribution function of ${\mu}$ (for more general $F_k$, a
slight modification of our random matrix model is enough to apply what
follows).

For each $N\geq 1$ and for each $k\geq 1$, let us consider a family
$M_1(N,k),\ldots, M_k(N,k)$ of independent replicas of the random matrix
$M$ of \eqref{18.06.08.1}, i.e.  $$M=U\left(
  \begin{array}{ccc}
    X_1&&\\
&\ddots&\\
&&X_N
  \end{array}
\right)U^*,
$$where $X_1,\ldots, X_N$ are independent random variables with cumulative distribution function  $F^\ff{k}$ and $U$ an unitary, Haar distributed,
$N\times N$
random matrix independent of $X_1,\ldots,X_ N$. 

\begin{itemize}
\item Let $k$ tend to infinity. Following the definition of
  ${\Lambda}_N^\vee(\mu)$ given in  theorem \ref{21.03.08.thm}, we get:
\be
\label{18.06.08.2} M_1(N,k)\vee\cdots\vee M_k(N,k)\kinf M_N \mbox{ in distribution,}
\ee 
where for all $N$, $M_ N$ is a   ${\Lambda}_N^\vee(\mu)$-distributed
random matrix. 

Let now 
$N$ tend to infinity. From theorem
\ref{bachelorette.17.04.08}, we deduce:
\be
\label{18.06.08.4}\operatorname{Empirical\; Spectral\; Law}(M_N)\dinf  {\Lambda}^\vee(\mu).
\ee
\item Now,
let  $N$ tend to infinity before $k$ does.  According to the law of large
numbers,  
the empirical
spectral law of  $M_i(N,k)$ converges almost surely, for each $k$ and 
$i=1,\ldots, k$, 
to the law with distribution function $F^\ff{k}$, which shall be denoted by
$\mu^\ff{k}$. Therefore, 
we deduce from 
the asymptotic freeness of independent randomly rotated random matrices (more precisely by the equality between formulas \eqref{15.06.08.3} and
\eqref{max.mat.15.06.08}), the following almost sure convergence: 
\be
\label{18.06.08.3}\operatorname{Empirical\; Spectral\; Law}(M_1(N,k)\vee\cdots\vee M_k(N,k))\dinf  \mu^\ff{k}\Box\hskip-9.3pt\vee\cdots \Box\hskip-9.3pt\vee\mu^\ff{k}.
\ee 
\end{itemize}

Joining \eqref {18.06.08.2},   \eqref {18.06.08.4},  \eqref
{18.06.08.3}, we get the following diagram 
$$
\begin{array}{ccc}M_1(d,k)\vee\cdots\vee M_k(d,k)&\; \;\;\stackrel{\;\;\quad\quad k\to \infty}{-\!\!\!-\!\!\!-\!\!\!-\!\!\!-\!\!\!\longrightarrow\!\!\!\!\!\!\!\!\!\!\!\!\!\!\!}& {\Lambda}_N^\vee(\mu)\\
\arrowvert&&\arrowvert\\
N\to\infty&&N\to\infty\\
\downarrow &&\downarrow\\ {\begin{array}{c}
\textrm{Empirical Spectral Law:}\\ \mu^\ff{k}\Box\hskip-9.3pt\vee\cdots \Box\hskip-9.3pt\vee\mu^\ff{k}\end{array}}& &{\begin{array}{c}
\textrm{Empirical Spectral Law:}\\
{\Lambda}^\vee(\mu)\end{array}}\end{array}
$$
The right-hand side of \eqref{sht.18.10.08.1}, i.e. the fact that $(\mu^\ff{k})^{\Box\hskip-5.6pt\vee k}\kinf \Lambda^\vee (\mu)$,  only means that one can add an edge  $\kinf$ between both bottom vertices of the diagram, i.e. that the operations $k\to \infty$ and $N\to \infty$ are commutative, which is quite expected.

\section{Proofs}

\subsection{Preliminary results}
Both of the theorems will be proved first for measures which are absolutely
continuous with respect to the Lebesgue measure, and then extended to all
measures by approximation. In this context, the appropriate  approximation
tool is given by the following lemma.
 
\begin{lem}\label{order.stats.15.04.08.2} Let us define, for $\mu$ \pro measure on the real line,  the function $F_\mu^{<-1>}$ on $(0,1)$ by $F_\mu^{<-1>}(u)=\min \{x\in \R \ste F_\mu(x)\geq u\}$.

(i) For $\mu,\nu$ \pro measures on the real line, we have the following
equalities (between quantities which can be infinite) $$\inf\{\eps >0\ste
F_\nu(\bullet-\eps)\leq F_\mu(\bullet)\leq
F_\nu(\bullet+\eps)\}=||F_\mu^{<-1>}-F_\nu^{<-1>}||_\infty=\inf
||X-Y||_\infty,$$ where the infimum in the third term is taken on pairs $(X,Y)$ of random variables defined on a same \pro space with respective distributions $\mu,\nu$.

(ii) For any \pro measure $\mu$ on the real line, for any $\eps >0$, there exists a \pro measure $\mu_\eps$ on  the real line   \st $F_{\mu_\eps}$ is smooth and $$F_{\mu_\eps}(\bullet-\eps)\leq F_\mu(\bullet)\leq F_{\mu_\eps}(\bullet+\eps).$$
\end{lem}
\begin{pr} (i) We prove these equalities by  cyclic majorizations.

- Let us prove first 
$$
\inf\{\eps >0\ste \forall x\in \R, F_\nu(x-\eps)\leq F_\mu(x)\leq
F_\nu(x+\eps)\}\geq ||F_\mu^{<-1>}-F_\nu^{<-1>}||_\infty.
$$
Let us consider $\eps>0$ \st $F_\nu(\bullet-\eps)\leq F_\mu(\bullet)\leq
F_\nu(\bullet+\eps)$. Since this is  equivalent to $F_\mu(\bullet-\eps)\leq
F_\nu(\bullet)\leq F_\mu(\bullet+\eps)$, this inequation is symmetric in
$\mu$ and $\nu$. Hence, it suffices to prove that
$F_\nu^{<-1>}(u)-F_\mu^{<-1>}(u)\leq\eps$ for all $u\in(0,1)$.
We have: 
$$
F_\nu(F_\mu^{<-1>}(u)+\eps)\geq F_\mu(F_\mu^{<-1>}(u))\geq u,
$$ 
hence, $F_\nu^{<-1>}(u)\leq F_\mu^{<-1>}(u)+\eps$. 

- The inequality 
$$
||F_\mu^{<-1>}-F_\nu^{<-1>}||_\infty\geq\inf ||X-Y||_\infty
$$ 
is due to the fact that for $U$ random variable with uniform distribution on $(0,1)$, $X:=F_\mu^{<-1>}(U), Y:=F_\nu^{<-1>}(U)$ are respectively distributed according to $\mu,\nu$. 

- To conclude, it suffices to prove 
$$
\inf ||X-Y||_\infty\geq \inf\{\eps >0\ste  F_\nu(\bullet-\eps)\leq
F_\mu(\bullet)\leq F_\nu(\bullet+\eps)\}.
$$
So let us consider a pair $(X,Y)$ of random variables defined on the same
\pro space with respective distributions $\mu, \nu$. Consider $\eps>0$ \st
$|X-Y|\leq \eps$ uniformly on the \pro space. Then, for all real number $x$, 
$$
\mathbb{P}(Y\leq x-\eps)\leq \mathbb{P}(X\leq x)\leq \mathbb{P}(Y\leq
x+\eps).
$$

(ii) Let ${\varepsilon}>0$, and $F$ a smooth non decreasing function on the real line   \st for all $k\in \z$, $F(k\eps)=F_\mu(k\eps)$. Such a function exists obviously, as an example, one can consider $F(x)=\int_{-\infty}^xf(t)\ud t$, for $$f(t)=\sum_{k\in\z}[F(k\eps)-F((k-1)\eps)]\vfi(t-k\eps)$$ with $\vfi$ smooth non negative function on the real line, with support contained in $[-\eps,0]$ satisfying  $\int_{-\eps}^0\vfi(t)\ud t=1$. Hence $F$ is a c\`adl\`ag non decreasing function on the real line \st $F(x)$ tends to zero (resp. one) as $x$ tends to $-\infty$ (resp. $+\infty$). It follows that $F=F_{\mu_\eps}$ for a certain \pro measure ${\mu_\eps}$ on the real line. Moreover, for all real number $x$, if $k\in \z$ is \st $k\eps\leq x <(k+1)\eps$, then, one has: 
$$
F_{\mu_\eps}(x-\eps)\leq F_{\mu_\eps}(k\eps)=F_\mu(k\eps)\leq F_\mu(x)\leq
F_\mu((k+1)\eps)=F_{\mu_\eps}((k+1)\eps)\leq F_{\mu_\eps}(x+\eps).
$$ 
\end{pr}

For any positive integer $ N$, let us define $\maxd$ to be the function from
$\cup_{n\geq N} \R^n$ to $\R^ N$ which maps any vector $x$ of $\R^n$, $n\geq
 N$, to the vector of the $ N$ largest coordinates of $x$ ranked in
decreasing order. The following property is a basic result,  see for
instance \cite{Da-Na}.

\begin{prop}\label{order.stats.15.04.08}
  Let $F$ be the cumulative distribution function of a \pro measure $\mu$ with no atom on the real line. Let, for $n\geq N$,  $X_1,\ldots, X_n$ be independent  random variables with law $\mu$. The  distribution of $\maxN(X_1,\ldots,X_n)$ has density 
$$(t_1,\ldots, t_N)\mapsto \f{n!}{(n-N)!}1_{t_1\geq\cdots\geq t_N}F(t_N)^{n-N}$$  
 with respect
to $\mu^{\otimes N}$.
\end{prop}

Before stating the next proposition, we shall recall that for $\mu$ be a \pro measure on the real line and $k$ positive integer, $\mu^\ff{k}$ is the \pro measure on the real line with cumulative distribution function $F_\mu^{\ff{k}}$, i.e. the law $m$ \st for $X_1,\ldots, X_k$ independent random variables with law $m$, $\max (X_1,\ldots, X_k)$ has law $\mu$.

\begin{prop}\label{JJ72.17.04.08}\textrm{(i)} Let $\mu$ be a \pro measure on the real line. For all integer $N\geq 1$, the push forward  of the \pro measure $ (\mu^\ff{k}
)^{\otimes kN}$ on $\R^{kN}$ by the function $\maxd$ converges weakly, as the integer $k$ tends to infinity,  to a \pro measure on $\R^N$ denoted by $\mu_ N$. 

\textrm{(ii)} When $\mu$ has no atom, $\mu_ N$ is absolutely continuous with
respect to  $\mu^{\otimes N}$, with density 
\begin{equation*}
N^N
1_{t_1\geq\cdots\geq
  t_N}\prod_{i=1} ^N\frac{F_\mu(t_N)}{F_\mu(t_i)}
\end{equation*}at any point $(t_1,\ldots, t_N)\in \R^ N$ \st for all
$i=1,\ldots,  N$, $F_\mu(t_i)\neq 0$.

\textrm{(iii)} Let us endow $\R^ N$ with the norm $||x||=\max_i |x_i|$. Then, for any \pro measures $\mu,\nu$ on the real line, $$\inf \{||V-W||_\infty\ste V,W \textrm{ random  vectors defined on the same  space with respective laws $\mu_N,\nu_ N$}\}$$ $$\leq  \inf \{||X-Y||_\infty\ste X,Y \textrm{ random  variables defined on the same  space with respective laws $\mu,\nu$}\}.$$

\textrm{(iv)} For any pair $\mu,\nu$ of \pro measures on the real line, if $\rho$ is the \pro measure on the real line \st $F_\rho=F_\mu F_\nu$, for all $N\geq 1$, $\rho_ N$ is the push-forward, by the function $\maxd$, of the \pro measure $\mu_N\otimes \nu_ N$ on $\R^{2N}$.
\end{prop}

\begin{pr} Let $k,N\geq 1$, and $\mu$ a probability measure with no atom.
  From property
  \ref{order.stats.15.04.08}, we infer first that $
  \mu^\ff{k}=\ff{k}F_\mu^{\ff{k}-1}\ud \mu$, and then
that the distribution $\maxN( (\mu^\ff{k}
)^{\otimes kN})$
 has a density with respect to $\mu^{\otimes N}$ equal to
$$(t_1,\ldots, t_N)\mapsto
\frac{(dk)!}{k ^N(dk-N)!}
1_{t_1\geq\cdots\geq
  t_N}\prod_{i=1} ^N\frac{F_\mu(t_N)^{1-\ff{k}}}{F_\mu(t_i)^{1-\ff{k}}}
$$at any point $(t_1,\ldots, t_N)\in \R^ N$ \st for all $i=1,\ldots,  N$,
$F(t_i)\neq 0$ (the density can obviously be set to zero anywhere else). 
As $k$ tends to infinity, this density stays uniformly bounded and
converges pointwise to 
\begin{equation*}
N^N
1_{t_1\geq\cdots\geq
  t_N}\prod_{i=1} ^N\frac{F_\mu(t_N)}{F_\mu(t_i)},
\end{equation*}
hence (i)  for \pro measures with no atom and (ii) are proved.

Now, let us complete the proof of  (i). Let us consider a \pro measure $\mu$ and, for all positive integer $k$, a family $X(k, 1), \ldots, X(k,kN)$ of independent random variables with law $\mu^\ff{k}$. To prove (i), by Theorem 1.12.4 of \cite{vanderVaart}, it suffices   to  prove that for any real bounded Lipschitz function $f$  on $\R^ N$, the sequence $\E(f(\maxN(X(k,1),\ldots, X(k,kN))))$ converges as $k$ tends to infinity. So let us fix such a function. We shall prove that the previous sequence is a Cauchy sequence. Let us fix $\eps>0$. 
Let us consider $\mu_\eps$ as in (ii) of lemma  \ref{order.stats.15.04.08.2}. 
 Note that for all positive integers $k$, we also have: 
$$
F_{\mu_\eps}^\ff{k}(\bullet-\eps)\leq F_\mu^\ff{k}(\bullet)\leq
F_{\mu_\eps}^\ff{k}(\bullet+\eps).
$$ 
So, based on (i) of lemma \ref{order.stats.15.04.08.2},   we may suppose that  for all $k\geq 1$, on the \pro space where the $X(k,i)$'s are defined, there is a family $$Y(k,1),\ldots,Y(k,kN)$$  of  random variables \st 

(a) the $Y(k,i)$'s are independent and distributed according to  $\mu_\eps^\ff{k}$,

(b) for all $i=1,\ldots, k N$, $|X(k,i)-Y(k,i)|\leq 2\eps$ almost surely. 

Note that $\mu_\eps$ has no atom, hence by (a) and what we just proved, the sequence $$\E(f(\maxN(Y(k,1),\ldots, Y(k,kN))))$$ is a Cauchy sequence. Note also that by (b), if $C$ is a Lipschitz constant for $f$ with respect to the norm $||x||=\max_i |x_i|$, then for all $k$, 
$$|\E(f(\maxN(X(k,1),\ldots, X(k,kN))))-\E(f(\maxN(Y(k,1),\ldots, Y(k,kN))))|\leq 2C\eps.$$ Hence there is $k_0\geq 1$ \st for all $k,k'\geq k_0$, 
$$|\E(f(\maxN(X(k,1),\ldots, X(k,kN))))-\E(f(\maxN(X(k',1),\ldots,
X(k',k'd))))|\leq (4C+1)\eps.$$Thus we have proved (i) for any probability
measure $\mu$. 

Now, let us prove (iii). Consider $\mu,\nu$ \pro measures on the real
line. Note that by part (i) of lemma \ref{order.stats.15.04.08.2}, it
suffices to prove that for all positive $\eps$ \st  $$
F_\nu(\bullet-\eps)\leq F_\mu(\bullet)\leq F_\nu(\bullet+\eps),$$ for all
$\alpha>\eps$,  there exists a pair $V,W$ of random vectors defined on the
same space with respective laws $\mu_N,\nu_ N$ \st $||V-W||_\infty\leq
\alpha$. Let us consider such a positive $\eps$ and $\alpha>\eps$. For all $k\geq
1$, we  also have: 
$$
F_{\nu}^\ff{k}(\bullet-\eps)\leq
F_\mu^\ff{k}(\bullet)\leq F_{\nu}^\ff{k}(\bullet+\eps).
$$ 
So, according to part (i) of lemma \ref{order.stats.15.04.08.2},   we shall consider, for all $k\geq 1$, a family 
$$X(k,1),\ldots, X(k,kN),Y(k,1),\ldots,Y(k,kN)$$  of  random variables defined on the same space \st 

(a) the $X(k,i)$'s are independent and distributed according to  $\mu^\ff{k}$,

(b) the $Y(k,i)$'s are independent and distributed according to  $\nu^\ff{k}$,

(c) for all $i=1,\ldots, k N$, $|X(k,i)-Y(k,i)|\leq \alpha$ almost surely. 

Let, for all $k\geq 1$, $\tau_k$ be the joint law, on $\R^{2N}$, of the random vector
$$(\maxN(X(k,1),\ldots, X(k,kN)),\maxN(Y(k,1),\ldots, Y(k,kN)).$$ The law
of $\maxN(X(k,1),\ldots, X(k,kN))$ (resp. of $\maxN(Y(k,1),\ldots,
Y(k,kN))$) converges weakly to $\mu_ N$ (resp. to $\nu_ N$) as $k$ tends to
infinity. It follows that the sequence $\tau_k$ is tight, i.e. relatively
compact for the topology of weak convergence (Theorem 6.1 of
\cite{MR0233396}). Let $\tau$ be the weak limit of a subsequence of
$\tau_k$. Let $(V,W)$ be a $\tau$-distributed random vector of $\R ^N\times
\R^ N$. Then the law of $V$ (resp. of $W$) is $\mu_ N$
(resp. $\nu_ N$). 
Moreover, it is easy to notice that for all $k$, $\tau_k$ is supported by $\{(x,y)\ste
x\in \R ^N,y\in \R ^N, ||x-y||\leq \alpha\}$. Indeed, let ${\sigma},{\tau}\in
S_n$ such that $X(k,{\sigma}(1))\geq \cdots\geq X(k,{\sigma}(kN))$ and
$Y(k,{\tau}(1))\geq \cdots\geq Y(k,{\tau}(kN))$; then for any
$i=1,\ldots,k N$, 
\beas
X(k,{\sigma}(i))&=&\max_{
    V\subset\left\{1,\ldots,kN\right\},
  \#V=i}\min(X(k,j),j\in V)\\ &\geq& \max_{V\subset\left\{1,\ldots,kN\right\},
  \#V=i}\min(Y(k,j)-{\alpha},j\in V)\\ &=&Y(k,{\tau}(i))-{\alpha}\quad\mbox{ a.s.}
\eeas
By symmetry, this proves $\left\vert
  X(k,{\sigma}(i))-Y(k,{\tau}(i))\right\vert\leq {\alpha}$ a.s., and this implies
\beas
\left\Vert\maxN(X(k,1),\ldots, X(k,kN))-\maxN(Y(k,1),\ldots, Y(k,kN)) \right\Vert&=&\max_{i=1,\ldots,N}\left\vert X(k,{\sigma}(i))-Y(k,{\tau}(i))\right\vert \\  &\leq&{\alpha} \quad \mbox{ a.s.}
\eeas
Letting $k$ go to infinity, this establishes the inequality  $\Vert V-W\Vert\leq \alpha$ almost surely.

To prove (iv), consider for any $k\geq 1$,
$$X(k,1),\ldots,X(k,kN),Y(k,1),\ldots,Y(k,kN)$$ is a  family of
independent random variables \st for all $i=1 ,\ldots, k N$, $X(k,i)$
(resp. $Y(k,i)$) is distributed according to $\mu^\ff{k}$
(resp. $\nu^\ff{k}$). Then $X(k,1)\vee Y(k,1),\ldots, X(k,kN)\vee Y(k,kN)$
are i.i.d. with law ${\rho}^{\frac{1}{k}}$, and consequently  $\mu_ N$,
$\nu_ N$ and $\rho_ N$ are the respective weak limit distributions, as $k$
tends to infinity, of 
$$
\maxd[X(k,1),\ldots, X(k,kN)],\quad \maxd[Y(k,1),\ldots, Y(k,kN)] 
$$  
\begin{equation}\label{21.03.08.2}\textrm{and }\quad \maxd[X(k,1)\vee Y(k,1),\ldots, X(k,kN)\vee Y(k,kN)].\end{equation}
Hence, by continuity and commutativity of the maximum operations,
$\maxN(\mu_N\otimes \nu_N)$ is the weak limit, as $k$ tends to infinity, of
the distribution of 
\begin{equation}\label{21.03.08.1}
\maxd[X(k,1),\ldots,X(k,kN),Y(k,1),\ldots,Y(k,kN)].
\end{equation}
 To conclude, we shall establish that   the distributions of the vectors
 of \eqref{21.03.08.2} and  \eqref{21.03.08.1} have the same weak limits as
 $k$ tends to infinity. Obviously, it  reduces to prove that the
 probability of the event  
$$
\hskip -4cm 
\{\maxd[X(k,1)\vee Y(k,1),\ldots, X(k,kN)\vee Y(k,kN)]
$$  
\begin{equation}\label{21.03.08.3}
\neq
  \maxd[X(k,1),\ldots,X(k,kN),Y(k,1),\ldots,Y(k,kN)]\}
\end{equation} 
tends to zero as $k$ tends to infinity. Note that this event is equivalent to the fact that there is $i\in \{1,\ldots, kN\}$ \st $X(k,i)$ and $Y(k,i)$ are both coordinates of $$\maxd[X(k,1),\ldots,X(k,kN),Y(k,1),\ldots,Y(k,kN)]$$ and are not equal to any of the others $X(k,j)$'s and $Y(k,l)$'s. Hence, if 
 one denotes by $I(k)$ (resp. $J(k)$) the set of $i$'s in $\{1,\ldots, kN\}$ \st $X(k, i)$ (resp. $Y(k,i)$) is one of the coordinates of $$\maxd[X(k,1),\ldots, X(k,kN)]\quad\textrm{ (resp. of $\maxd[Y(k,1),\ldots, Y(k,kN)]$),}$$  the event of \eqref{21.03.08.3}
 is contained in the event $I(k)\cap J(k)\neq \emptyset$, which \pro is equal to ${\binom{kN}{N}}^{-1}$,   
since we may choose for $I(k)$ and $J(k)$ to be independent and to have uniform distribution on the set of subsets of $\{1,\ldots, kN\}$ with cardinality $ N$. Thus, (iv) is proved.\end{pr}

\subsection{Proof of proposition \ref{21.03.08.propo1}}The first statement is immediate. The third one follows
  immediately form the second one. To establish the second statement, it
  suffices to prove  that  almost surely, for any $t\in\R$,
  $\im(\chi_{(t,\infty)}(A))+\im(\chi_{(t,\infty)}(B))$ has the maximal
  dimension conditionally to the dimensions of $\im(\chi_{(t,\infty)}(A)$
  and $\im(\chi_{(t,\infty)}(B)$,  i.e. that 
  \begin{equation}
    \label{eq:dpq}
\operatorname{dim}[\im(\chi_{(t,\infty)}(A))+\im(\chi_{(t,\infty)}(B)) ]=\min \{N,
\operatorname{dim}[\im(\chi_{(t,\infty)}(A))]+\operatorname{dim}[\im(\chi_{(t,\infty)}(B)]\}\quad\mbox{ a.s.}    
  \end{equation}

Let $p=\operatorname{dim}[\im(\chi_{(t,\infty)}(A))]$,
$q=\operatorname{dim}[\im(\chi_{(t,\infty)}(B))]$. Due to the unitarily
invariance of $A$ (resp. $B$), the law of
$\im(\chi_{(t,\infty)}(A))$
(resp. $\im(\chi_{(t,\infty)}(B))$) is also invariant under the
action of any unitary matrix, hence is uniform on the set of all subspaces
of dimension $p$ (resp. $q$). Moreover,  $\im(\chi_{(t,\infty)}(A))$ and
$\im(\chi_{(t,\infty)}(B))$ are independent.

Let us consider a family $(g_{i,j},i\in\left\{1,\ldots,N\right\},j\geq 1)$ of independent complex standard Gaussian
random variables. The matrix $(g_{i,j},i,j\in\left\{1,\ldots,N\right\})$ is
distributed w.r.t. the Circular Unitary Ensemble. Hence, it is unitarily invariant
and a.s. inversible. Therefore, the subspace generated by 
vectors $(g_{i,j},i\in\left\{1,\ldots,N\right\})_{j=1}^p$
(resp. $(g_{i,j},i\in\left\{1,\ldots,N\right\})_{j=p+1}^{p+q}$, $(g_{i,j},i\in\left\{1,\ldots,N\right\})_{j=1}^{p+q})$ has dimension
$p$ (resp. $q$, $\min(N,p+q)$) a.s.,
and its law is uniform on the set of all subspaces
of same dimension. This proves $(\ref{eq:dpq})$.

\begin{rmq}Note that the proof of the second statement could also have been deduced from Theorem 2.2 of \cite{MR2198015}.\end{rmq}

\subsection{Proof of theorem \ref{21.03.08.thm}}Following proposition \ref{21.03.08.propo1},  if we denote $M^{(k)}=M_1\vee\cdots\vee M_k$, then $M^{(k)}$ is
an Hermitian random matrix, whose law is invariant under conjugation with any
unitary matrix and whose eigenvalues ${\lambda}^{(k)}_1\geq\ldots\geq {\lambda}^{(k)}_ N$ are  equal
to the $ N$ largest of the $Nk$ eigenvalues of $M_1,\ldots,M_k$ which are all
independent with same distribution   $\mu^{1/k}$. Therefore, by theorem 4.3.5 of \cite{H-P},
  there exists $U$ unitary, Haar distributed and independent of the vector $( {\lambda}^{(k)}_1,\ldots, {\lambda}^{(k)}_N)$, such that:  
$$
M^{(k)}=U\left(
  \begin{array}{ccc}
    {\lambda}^{(k)}_1&&\\
&\ddots&\\
&&{\lambda}^{(k)}_N
  \end{array}
\right)U^*.
$$
Note that the distribution of the vector $( {\lambda}^{(k)}_1,\ldots, {\lambda}^{(k)}_N)$ is  the push forward  of the \pro measure $ (\mu^\ff{k}
)^{\otimes kN}$ on $\R^{kN}$ by the function $\maxd$. Hence, the first part of the theorem follows from (i) and (ii) of proposition \ref{JJ72.17.04.08},   and the second part follows from (iv) of the same proposition.

\subsection{Proof of corollary \ref{ultima}}
Both assertions are easy consequences of the last part of theorem \ref{21.03.08.thm}.
\begin{itemize}
\item[(i)] Let $p\geq 1$. By definition, 
  $F_{\mu}=\left(F_{{\mu}^{\frac{1}{p}}}\right)^p$. Therefore
    ${\Lambda}^\vee_N(F_{\mu})$ is the distribution of $H_1\vee\cdots\vee
    H_p$, with $H_1,\ldots,H_p$ i.i.d. with law
    ${\Lambda}^\vee_N(F_{\mu^{\frac{1}{p}}})$. This means that ${\Lambda}^\vee_N(F_{\mu})$ is max-infinitely divisible.
\item[(ii)]  Let $p\geq 1$. Since ${\mu}$ is max-stable, there exists $a_p>0$,
  $b_p\in\R$ such that ${\mu}$ is the distribution of
  $(a_pX_1+b_p)\vee\cdots\vee(a_pX_p+b_p)$, where $X_1,\ldots,X_p$
  are i.i.d. with law ${\mu}$. Let $H_1,\ldots,H_p$ be i.i.d. with law
  ${\Lambda}^\vee_N({\mu})$. It is quite immediate from the construction
  that the image of the distribution of $a_pX_i+b_p$ is the law of
  $a_pH_i+b_bI_ N$ where $I_ N$ denotes the identity matrix. Therefore,
  ${\Lambda}^\vee_N({\mu})$ is the distribution of
  $(a_pH_1+b_pI_N)\vee\cdots\vee(a_pH_p+b_pI_N)$. This proves that
  ${\Lambda}^\vee_N({\mu})$ is max-stable.
\end{itemize}

\subsection{Proof of theorem \ref{bachelorette.17.04.08}}
Suppose first that $F_\mu$ is smooth on the real line. It implies that $\mu$ is absolutely continuous with respect to the Lebesgue measure and that $\ud \mu(x)=F_{\mu}'(x)\ud x$.
 Since $F_\mu$ is continuous on the real line, the almost sure weak convergence of the sequence of random \pro measures $\hat\mu_ N$ to $\Lambda^\vee(\mu)$ is equivalent to the fact the for all rational number $t$,  $\hat{\mu}_N\left((-\infty,t]\right)$
converges almost surely to $\max(0,1+\log F_\mu(t))$. We shall prove that for all $t\in\R$, $\hat{\mu}_N\left((-\infty,t]\right)$
converges almost surely to $\max(0,1+\log F_\mu(t))$. Let us fix a real number $t$.

The first step is to
compute
$  \E\left[e^{{\lambda}\hat{\mu}_N\left((-\infty,t]\right)}\right]$: with the notations ${\lambda}_0=+\infty$, ${\lambda}_{N+1}=-\infty$, we have
\begin{eqnarray*}
  \E\left[e^{{\lambda}\hat{\mu}_N\left((-\infty,t]\right)}\right]&=&  \E\left[\sum_{p=0} ^Ne^{\frac{{\lambda}(N-p)}{N}} 1_{{\lambda}_{p+1}\leq t<{\lambda}_p}\right]\\ &=&\sum_{p=0}^Ne^{\frac{{\lambda}(N-p)}{N}}\P\left({\lambda}_{p+1}\leq t<{\lambda}_p\right)
\end{eqnarray*}
Now, note that for all real number $T$ \st $F(T)>0$, for all integer non negative $m$, $$\int_{(T, \infty)}(\log F_\mu(t))^m\frac{\ud \mu(t)}{F_\mu(t)}=-\frac{(\log F_\mu(T))^{m+1}}{m+1}\textrm{ and }\quad\int_{(-\infty, T]} F_\mu(t)^m{\ud \mu(t)}=\frac{ F_\mu(T)^{m+1}}{m+1}.$$ Hence for
$p=0,\ldots,N-1$,

$$ \P\left({\lambda}_{p+1}\leq t<{\lambda}_p\right)=\int_{t\in\R ^N}N^N1_{t\geq
    t_{p+1}\geq\cdots\geq t_N}1_{t<t_p\leq\cdots\leq
    t_1}\prod_{i=1}^N\frac{ F_\mu(t_N)}{F_\mu(t_i)}\ud\mu(t_i)=
\frac{N^p}{p!}F_\mu(t)^N(-\log F_\mu(t))^p.$$
And for $p= N$:
$$
  \P\left(t<{\lambda}_{N}\right)=\frac{N^N}{(N-1)!}\int_{t}^{+\infty}F_\mu(t_N)^N\left(-\log F_\mu(t_N)\right)^{N-1}\ud \mu(t_N)
=\int_{0}^{-\log F_\mu(t)}\frac{N^N}{(N-1)!}u^{N-1}e^{-Nu}\ud u.$$
Thus, using the fact that
$\sum_{p=0}^{N-1}\frac{x^p}{p!}=e^x\int_{x}^{+\infty}\frac{1}{(N-1)!}u^{N-1}e^{-u}\ud u$,
we get:
\begin{eqnarray*}
  &&\hskip-1cm  \E\left[e^{{\lambda}\hat{\mu}_N\left((-\infty,t]\right)}\right]\\
&=&e^{\lambda}F_\mu(t)^N\sum_{p=0}^{N-1}\frac{(-Ne^{-\frac{\lambda}{N}}\log
  F_\mu(t))^p}{p!}+\int_{0}^{-\log F_\mu(t)}\frac{N^N}{(N-1)!}u^{N-1}e^{-Nu}\ud u\\
&=&e^{\lambda}F_\mu(t)^Ne^{-Ne^{-\frac{\lambda}{N}}\log
  F_\mu(t)}\int_{-Ne^{-\frac{\lambda}{N}}\log
  F_\mu(t)}^{+\infty}\frac{1}{(N-1)!}u^{N-1}e^{-u}\ud u+\int_{0}^{-\log
  F_\mu(t)}\!\!\!\!\!\!\!\!\!\frac{N^N}{(N-1)!}u^{N-1}e^{-Nu}\ud u\\
&=&e^{{\lambda}+N\left(1-e^{-\frac{\lambda}{N}}\right)\log
  F_\mu(t)}\int_{-e^{-\frac{\lambda}{N}}\log
  F_\mu(t)}^{+\infty}\frac{N^N}{(N-1)!}u^{N-1}e^{-Nu}\ud u+\int_{0}^{-\log
  F_\mu(t)}\frac{N^N}{(N-1)!}u^{N-1}e^{-Nu}\ud u.
\end{eqnarray*}
The following inequalities will be useful, and they can be established from the  Chernov inequality or proved directly:
\begin{equation}\label{rock.with.you.17.03.08}
\forall x\in[0,1),  \int_{0}^{x}\frac{N^N}{(N-1)!}u^{N-1}e^{-Nu}\ud u\leq x^Ne^{N(1-x)},\quad \forall x>1, 
  \int_{x}^{+\infty}\frac{N^N}{(N-1)!}u^{N-1}e^{-Nu}\ud u\leq x^Ne^{N(1-x)}.\end{equation}
We have three cases to consider:

{\bf Case 1:}  $-\log F_\mu(t)>1$. If moreover
  ${\lambda}\in[0,N\log\left(-\log F_\mu(t)\right))$, then $-e^{-\frac{\lambda}{N}}\log
  F_{\mu}(t)>1$, and by \eqref{rock.with.you.17.03.08},
  \begin{eqnarray*}
&&e^{{\lambda}+N\left(1-e^{-\frac{\lambda}{N}}\right)\log
  F_\mu(t)}\int_{-e^{-\frac{\lambda}{N}}\log
  F_\mu(t)}^{+\infty}\frac{N^N}{(N-1)!}u^{N-1}e^{-Nu}\ud u\\
&\leq&e^{{\lambda}+N\left(1-e^{-\frac{\lambda}{N}}\right)\log
  F_\mu(t)}\left(-e^{-\frac{\lambda}{N}}\log
  F_\mu(t)\right)^Ne^{N\left(1+e^{-\frac{\lambda}{N}}\log
  F_\mu(t)\right)}\\
&=&e^{N\left(1+\log F_\mu(t)+\log(-\log F_\mu(t))\right)}<1
  \end{eqnarray*}
since $1+\log F_\mu(t)+\log(-\log F_\mu(t))<0$.
Hence for all $N\geq 1$,  $  \E\left[e^{{\lambda}\hat{\mu}_N\left((-\infty,t]\right)}\right]\leq
2$,
and from the Chernov
inequality, we get then: for all ${\varepsilon}>0$,
$$
\P\left\{\hat{\mu}_N((-\infty,t])>{\varepsilon}\right\}\leq\inf_{{\lambda}\in[0,N\log\left(-\log F_\mu(t)\right))}2e^{-{\varepsilon}{\lambda}}=2e^{-{\varepsilon}N\log\left(-\log F_\mu(t)\right))}
$$
This is enough, with Borel-Cantelli lemma, to prove the almost-sure
convergence of $\hat{\mu}_N((-\infty,t])$ to $0=\max(0,1+\log F(t))$.

{\bf Case 2:}   $-\log F_\mu(t)<1$. Then by \eqref{rock.with.you.17.03.08},
$$
\int_{0}^{-\log
  F_\mu(t)}\frac{N^N}{(N-1)!}u^{N-1}e^{-Nu}\ud u\leq \left(-\log
  F_\mu(t)\right)^Ne^{N\left(1+\log F_\mu(t)\right)}=e^{N{\delta}_t}
$$
with ${\delta}_t= 1+\log F_\mu(t)+\log(-\log F_\mu(t))<0$. 
\begin{itemize}
\item Moreover, we suppose  ${\lambda}\geq 0$. Then
\begin{eqnarray*}
  \E\left[e^{
  {\lambda}\left[\hat{\mu}_N\left((-\infty,t]\right)-(1+\log
      F_\mu(t))\right]
     }
      \right]
&\leq&\left(e^{\lambda}F_\mu(t)^Ne^{-Ne^{-\frac{\lambda}{N}}\log
  F_\mu(t)}+e^{N{\delta}_t}\right)e^{-{\lambda}(1+\log F_\mu(t))}\\
&=&e^{-N\log
  F_\mu(t)\left(\frac{\lambda}{N}-(1-e^{-\frac{{\lambda}}{N}})\right)}+e^{N\left({\delta}_t-\frac{\lambda}{N}(1+\log
  F_\mu(t))\right)}
\end{eqnarray*}
From the Chernov inequality, we infer then: for all ${\varepsilon}>0$
\begin{eqnarray*}
\P\left\{\hat{\mu}_N((-\infty,t])-(1+\log
    F_\mu(t))>{\varepsilon}\right\}
&\leq&\inf_{{\lambda}\geq 0}\left(e^{-N\log
  F_\mu(t)\left(\frac{\lambda}{N}-(1-e^{-\frac{{\lambda}}{N}})\right)}+e^{N\left({\delta}_t-\frac{\lambda}{N}(1+\log
  F_\mu(t))\right)}\right)e^{-{\lambda}{\varepsilon}}\\ 
&=&\inf_{{\lambda}\geq 0}\left(e^{-N\left(\log
  F_\mu(t)\left({\lambda}-(1-e^{-{\lambda}})\right)+{\lambda}{\varepsilon}\right)}+e^{N\left({\delta}_t-{\lambda}(1+\log
  F_\mu(t)+{\varepsilon})\right)}\right)\\
&\leq&2e^{-NC}
\end{eqnarray*}
for some constant $C>0$.
This gives, using Borel-Cantelli lemma:
$$
\limsup_{N\fleche+\infty}\hat{\mu}_N((-\infty,t])\leq 1+\log
    F_\mu(t)\mbox{ a.s.}
$$
\item We suppose now ${\lambda}\leq 0$. Using the same trick, we get
  \begin{eqnarray*}
\P\left\{\hat{\mu}_N((-\infty,t])-(1+\log
    F_\mu(t))<-{\varepsilon}\right\}
&\leq&\inf_{{\lambda}\leq 0}\left(e^{-N\left(\log
  F_\mu(t)\left({\lambda}-(1-e^{-{\lambda}})\right)-{\lambda}{\varepsilon}\right)}+e^{N\left({\delta}_t-{\lambda}(1+\log
  F_\mu(t)-{\varepsilon})\right)}\right)\\
&\leq&2e^{-NC'}    
  \end{eqnarray*}
for some constant $C'>0$. This gives:
$$
\liminf_{N\fleche+\infty}\hat{\mu}_N((-\infty,t])\geq 1+\log
    F_\mu(t)\quad \mbox{ a.s.}
$$
\end{itemize}
And the result follows:
$$
\lim_{N\fleche+\infty}\hat{\mu}_N((-\infty,t])= 1+\log
    F_\mu(t)=\max(0,1+\log
    F_\mu(t))\quad\mbox{ a.s.}
$$

{\bf Case 3:}  $\log F_\mu(t)=-1$. We know that $F_\mu(u)$ tends to $0$ as $u$ tends to $-\infty$. Hence, there is $u<t$ \st $\log F_\mu(u)<-1$. For all $ N$,  $\hat{\mu}_N((-\infty,t])\geq \hat{\mu}_N((-\infty,u])$, hence, $$
\liminf_{N\fleche+\infty}\hat{\mu}_N((-\infty,t])\geq \liminf_{N\fleche+\infty}\hat{\mu}_N((-\infty,u])=0=1+\log F_{\mu}(t)\quad\mbox{ a.s.}$$
Moreover, for any positive $\varepsilon$, by  continuity  of $F_\mu$ and by the fact that $F_\mu(u)$ tends to $1$ as $u$ tends to $+\infty$, there is $u>t$ such that $1+\log F_\mu(t)<1+\log F_\mu(u)<1+\log F_\mu(t)+\varepsilon.$
  Hence, 
$$
\limsup_{N\fleche+\infty}\hat{\mu}_N((-\infty,t])\leq\limsup_{N\fleche+\infty}\hat{\mu}_N((-\infty,u])=1+\log
    F_\mu(u) <1+\log F_\mu(t)+\varepsilon\mbox{ a.s.}
$$
Therfore, $\lim_{N\fleche+\infty}\hat{\mu}_N((-\infty,t])=1+\log F_\mu(t)$
almost surely.

The theorem is proved in the case where $F_\mu$ is smooth. Now, we consider a \pro measure $\mu$ on the real line without making any assumption about $F_\mu$. By Theorem 1.12.4 of \cite{vanderVaart}, the distance $$\ud\,(m,m'):=\sup\left|\int f\ud m-\int f\ud m' \right|,$$where the $\sup$ is taken on the set $BL_1$
bounded Lipschitz functions $f$ on the real line with Lipschitz constant $\leq 1$ and \st $||f||_\infty\leq 1$, is a distance which defines the weak topology on the set of \pro measures on the real line.

Thus we have to prove that almost surely, as $ N$ tends to infinity,
$$
\lim\ud\,(\hat{\mu}_N,\Lambda^\vee({\mu}))=0,
$$
i.e. that for all $\eps>0$, almost surely, for $ N$ large enough,  
\begin{equation}\label{18.04.08.1} 
\quad\ud\,(\hat{\mu}_N,\Lambda^\vee(\mu))\leq \eps.
\end{equation}
So let us fix $\eps>0$. By part (ii) of lemma \ref{order.stats.15.04.08.2},
there exists a \pro measure $\nu$ on  the real line   \st $F_{\nu}$ is
smooth and 
$$
F_{\nu}(\bullet-\eps/6)\leq F_\mu(\bullet)\leq F_{\nu}(\bullet+\eps/6).
$$
Note that the same obviously holds if one replaces $\mu$ by
$\Lambda^\vee(\mu)$ and $\nu $ by $\Lambda^\vee(\nu)$, since for all \pro
measure $m$,  $F_{\Lambda^\vee(m)}=\max(0,1+\log F_m)$. Hence, according to
 part (i) of lemma \ref{order.stats.15.04.08.2}, 
\begin{equation}\label{anotherpartofme.18.04.08} 
\ud\,(\Lambda^\vee(\mu), \Lambda^\vee(\nu))\leq \eps/6.
\end{equation}
Note also that by definition of $\Lambda_N^\vee(\mu)$,    the vector
$(X_1\geq\cdots\geq X_N)$ of  ranked eigenvalues of $M_ N$ has law
$\mu_ N$. Following part (iii) of proposition \ref{JJ72.17.04.08}, one can suppose that for all $ N$, on the same space as $M_ N$, there is a random matrix $N_ N$ with law $\Lambda_N^\vee(\nu)$ with ranked eigenvalues $(Y_1\geq \cdots\geq Y_k)$ \st for all $i=1,\ldots,  N$, $|X_i-Y_i|\leq \eps/3$ almost surely. 
 If one denotes  the spectral law of $N_ N$ by $\hat{\nu}_ N$, it implies that
 \begin{equation}\label{anotherpartofme.18.04.08.2} \forall N\geq  1, \ud\,(\hat{\mu}_N,\hat{\nu}_N)\leq \eps/3\quad\textrm{ almost surely.}\end{equation}
Please note finally that by the first part of the proof,  \begin{equation}\label{anotherpartofme.18.04.08.3} \textrm{almost surely, for $ N$ large enough, }\quad\ud\,(\hat{\nu}_N,\Lambda^\vee(\nu))\leq \eps/3.\end{equation}
 Equations \eqref{anotherpartofme.18.04.08}, \eqref{anotherpartofme.18.04.08.2}, \eqref{anotherpartofme.18.04.08.3},  together, imply \eqref{18.04.08.1}. Thus, the theorem is proved.

\subsection{Proof of the implication \eqref{sht.18.10.08.1}}
Let us end the paper with the proof of the implication \eqref{sht.18.10.08.1}. We consider a sequence $F_k$ of cumulative distribution functions and a cumulative distribution function $F$ \st at any $x$ where $F$ is continuous, $F_k(x)^k\kinf F(x)$. Let us prove that for any such $x$, $ F_k(x)^{\Box\hskip-5.6pt\vee k}\kinf \Lambda^\vee (F(x)).$ Let us denote, for $y\in \R$, $y^+=\max(0,y)$. Note that since for any $a,b\in \R$ \st $b\leq 0$, we have $(a^++b)^+=(a+b)^+$, by induction on $k$, we prove easily that for all $G_1,\ldots, G_k$ distribution functions, $$G_1\Box\hskip-9.3pt\vee\cdots 
\Box\hskip-9.3pt\vee G_k=(G_1+\cdots +G_k-k+1)^+.$$Thus we have to prove $$(1+k(F_k(x)-1))^+\kinf (1+\log F(x))^+,\quad\textrm{  i.e. }\quad k(F_k(x)-1)\kinf \log F(x),$$ which follows directly from the hypothesis $F_k(x)^k\kinf F(x)$.
\bibliographystyle{plain}
\bibliography{BP-max}

\end{document}